\documentclass[12pt]{amsart}

\usepackage{amsmath,amsfonts,amssymb,amscd}
\begin{document}
\renewcommand{\thefootnote}{\fnsymbol{footnote}}
\pagestyle{plain}

\title{Chow-stability 
and Hilbert-stability\\
in Mumford's\\
Geometric Invariant Theory}
\author{Toshiki Mabuchi${}^*$}
\maketitle
\footnotetext{ ${}^*$
Special thanks are due to Professors Akira Fujiki and Sean T. Paul for 
useful comments during the preparation of this paper. }
\abstract
In this note, we shall show that 
Chow-stability and Hilbert-stability in GIT asymptotically coincide. 
The proof in \cite{Ma0} is simplified in the present form, 
while a quick review 
is in \cite{Mab}.
\endabstract
\section{Introduction}
For moduli spaces of polarized algebraic varieties, 
a couple of stability concepts are known in algebraic geometry (cf.~Mumford et al.~\cite{Mu}):
Chow-stability and Hilbert-stability.
In this note, we clarify the asymptotic relationship between 
them.
Throughout this note, we fix once for all a very ample holomorphic line bundle $L$ over 
an irreducible projective algebraic variety $M$ defined over $\Bbb C$. Let $n:= \dim M > 0$
and let $\ell$ be a positive integer with $\ell \geq n+1$.
Replacing $L$ by its suitable power, we may assume that
$H^i(M, \mathcal{O}(L^j)) \, =\,\{0\}$
for all positive integers $i$ and $j$.
Then associated to the complete linear system $|L^\ell|$,
we have the Kodaira embedding 
$$
\iota_\ell  : M\; \hookrightarrow \;\Bbb P^*(V_\ell),
$$
where $\Bbb P^*(V_\ell)$ is the set of all hyperplanes in 
$V_\ell := H^0(M, \mathcal{O}(L^\ell))$ through the origin.
Let $n$ and $d_\ell$ be respectively the dimension of $M$ and the
degree of  $\iota_\ell
(M )$ in the projective space $\Bbb P^* (V_\ell)$. 
Put $G_\ell := \operatorname{SL}_{\Bbb C}(V_\ell)$ 
and $W_\ell :=\{S^{d_\ell} (V_\ell )\}^{\otimes n+1}$, 
where $S^{d_\ell} (V_\ell )$ 
denotes the 
$d_\ell$-th symmetric tensor product of the space $V_\ell$.
Take an element
$M_\ell\neq 0$ in $W_\ell^*$ such that 
the associated element $[M_\ell ]$ in 
$\Bbb P^* (W_{\ell})$ is the Chow point of the 
irreducible reduced algebraic cycle
$\iota_\ell (M)$ on $\Bbb P^* (V_\ell )$. 
For the natural action of $G_{\ell}$ on $W_\ell^*$,
let $\hat{G}_{\ell}$ denote the isotropy subgroup of $G_{\ell}$ 
at $M_\ell$.

\medskip\noindent
{\em Definition $1.1$}.
(a)\, $(M,L^\ell)$ is called {\it Chow-stable\/} or {\it Chow-semistable}
according as the orbit $G_\ell\cdot M_\ell$ is closed in $W_\ell^*$ 
with $|\hat{G}_{\ell}| <\infty$
or the closure of $G_\ell\cdot M_\ell$  in $W_\ell^*$ is disjoint from the origin.
\newline
(b)\, $(M, L)$ is called {\it asymptotically Chow-stable\/} if 
$(M, L^\ell)$ is Chow-stable for all $\ell \gg 1$.

\medskip
Let $\ell$ and $k$ be positive integers.
Then the kernel $I_{\ell,k}$ of the natural homomorphism 
of $S^k (V_\ell )$ to $V_{\ell k} := H^0(M, \mathcal{O}_M(L^{\ell k}))$
is the degree $k$ component of
the homogeneous ideal defining $M$ in 
$\Bbb P^*(V_\ell)$. Put $m_k := \dim V_{\ell k}$ and $\gamma_{\ell,k} := \dim I_{\ell,k}$.
Then $\wedge^{\gamma_{\ell,k}} I_{\ell,k}$ is a complex line 
in $F_{\ell,k}:=\wedge^{\gamma_{\ell,k}} (S^k (V_\ell ))$.
Take an element $f_{\ell,k}\neq 0$ in $\wedge^{\gamma_{\ell,k}} I_{\ell,k}$. 
For the natural action of $G_{\ell}$ on $F_{\ell, k}$, let $\hat{G}_{\ell, k}$ be the 
isotropy subgroup of $G_{\ell}$ at $f_{\ell,k}$. 

\medskip\noindent
{\em Definition $1.2$}.  (a)\, $(M,L^\ell)$ is called {\it Hilbert-stable} if
the orbit $G_\ell \cdot f_{\ell,k}$ is closed in $F_{\ell,k}$
with $|\hat{G}_{\ell, k}| <\infty$ for all $k \gg 1$.
\newline
 (b)\, $(M,L)$ is called {\it asymptotically 
Hilbert-stable} if $(M,L^{\ell})$ is Hilbert stable
for all $\ell \gg 1$.

\medskip
A result of Fogarty \cite{Fo} (see also \cite{Mu}, p.215) states that
Chow-stability for $(M, L^\ell)$ implies  
Hilbert-stability for $(M, L^\ell)$. 
However, little was known for the converse implication.

\medskip
Consider the maximal connected linear algebraic subgroup $H$ 
of the group of holomorphic automorphisms of $M$.
To each positive integral multiple $L^{m}$ of $L$, we associate 
the point $[L^{m}] \in\operatorname{Pic}(M)$ 
defined by
$L^{m}$.  
For the natural $H$-action on $\operatorname{Pic}(M)$,
we denote by $\hat{H}_{m}$ the identity component of 
the isotropy subgroup
 of $H$ at $[L^{m}]$. Put $\hat{H}:= \hat{H}_1$.
Since the orbit $\hat{H}_m \cdot [L]\; (\cong \hat{H}_m/\hat{H})$
sitting in $\{L'\in \operatorname{Pic}(M)\,;\, (L')^m = L^m \}$ reduces to a single point, 
we have
$$
\hat{H}\; =\; \hat{H}_m  \qquad\text{for all $m\in \Bbb Z_+$.}
$$
Let $\{\,k_i\,; i =0,1,2, \dots\,\}$ be a sequence  of 
integers $\geq n+1$.
For a positive integer $\ell$, 
we define a sequence  
$\{\ell_i\}$ of positive integers
inductively by  setting
$\ell_{i+1} := \ell_i\, k_i$ and $\ell_0 := \ell$.
In this paper, we shall show that 

\medskip\noindent
{\bf Main Theorem.} 
{\em
{\rm (a)}\, 
Assume that $G_{\ell_i} \cdot f_{\ell_i ,k_i}$ is closed in $F_{\ell_i ,k_i}$ for all 
integers $i \geq 0$. If $\hat{H} = \{1\}$,
then $G_{\ell}\cdot M_{\ell}$ is closed in $W_{\ell}^*$.
\newline
{\rm (b)}\, $(M,L)$ is asymptotically Chow-stable if and only if 
$(M,L)$ is asymptotically Hilbert-stable. 
}

\medskip
As seen in the beginning of Section 3, (b) follows from (a). Hence, we here sketch the proof 
of (a) of Main Theorem. Assume $\hat{H} = \{1\}$.
Since $G_{\ell_i}\cdot f_{\ell_i,k_i}$ is closed in $F_{\ell_i,k_i}$ for all $i$, 
Lemma 3.10 shows that the polynomial Hilbert weight $w_{\lambda}= w_{\lambda}(k;\ell)$ 
in Section 3 is increasing
$$
0<w_{\lambda}(K_0; \ell )<w_{\lambda}(K_1; \ell )< \dots
< w_{\lambda}(K_{i-1} ; \ell )<  w_{\lambda}(K_i ; \ell )
< \dots
$$
for $K_i$, $i$=1,2,\dots, in (3.7), 
where $\lambda : \Bbb C^* \hookrightarrow G_{\ell}$ is an arbitrary 
algebraic one-parameter subgroup.
Since the asymptotic limit 
$$
w_{\lambda}(\infty ;\ell ):= \lim_{k\to\infty} w_{\lambda}(k;\ell )
$$
always exist, and since $K_i \to +\infty$ as $i \to \infty$, 
we have $w_{\lambda}(\infty ;\ell ) > 0$.  
This means that $(M,L^{\ell})$ is Chow-stable, i.e., $G_{\ell}\cdot M_{\ell}$ 
is closed in $W_{\ell}$.

\medskip
This paper is organized as follows.
First, Section 2 is given as a preparation for Section 5. 
Then the proof of Main Theorem will be outlined in Section 3, while 
two main difficulties (3.6) and Lemma 3.10
will be treated in Sections 4 and 5, respectively.
Finally, (5.1) is a key in the proof of Lemma 3.10,
and will be discussed in Appendix.

\medskip\noindent
{\em Acknowledgements}: I here thank the referee for his careful reading of the paper and for his numerous
suggestions.


\medskip
\section{A test configuration and the group action $\rho_k$}

Hereafter, 
we fix an action of an algebraic torus $T := \Bbb C^*$ 
on $\Bbb A^1 := \{\,s \,;\, s\in \Bbb C\,\}$ by multiplication 
of complex numbers
$$
T\times \Bbb A^1 \to \Bbb A^1, 
\qquad (t, s) \mapsto
ts.
$$
Let $\pi : \mathcal{Z} \to \Bbb A^1$ be a $T$-equivariant projective morphism
between complex varieties with a relatively very ample invertible sheaf
$\mathcal{L}$ on $\mathcal{Z}$ over $\Bbb A^1$, where
the algebraic group $T$ acts on $\mathcal{L}$, linearly on fibers, lifting
the $T$-action on $\mathcal{Z}$. 
Now the following concept by Donaldson will play a very important role 
in our study:

\medskip\noindent
{\em Definition $2.1$} (cf. \cite{D1}).  
\, $\pi : \mathcal{Z} \to \Bbb A^1$ above
is called a {\it test configuration\/} of exponent $\ell$ for $(M,L)$ if,
when restricted to fibers $\mathcal{Z}_s := 
\pi^{-1}(s)$, we have isomorphisms
$$
(\mathcal{Z}_s, {\mathcal{L}}_{|\mathcal{Z}_s} )
\cong (M, \mathcal{O}_M( L^{\otimes \ell})),
\qquad 0 \neq s \in \Bbb A^1.
$$

\medskip
Let $\pi : \mathcal{Z} \to \Bbb A^1$ be a test configuration of exponent $\ell$ 
for $(M,L)$. To each positive integer $k$,
we assign a vector bundle $E_k$ over $\Bbb A^1$ 
associated to the locally free sheaf
$\pi_*\mathcal{L}^k$ over $\Bbb A^1$, i.e.,  
$\mathcal{O}_{\Bbb A^1}(E_k) = \pi_*\mathcal{L}^k$.
For the natural $T$-action
$$
\rho_k : T \times {E}_k \to {E}_k
$$
induced by the $T$-action on $\mathcal{L}$,
we denote by $\rho_{k,0}$ the restriction
of the $T$-action
$\rho_k$ to the fiber $({E}_k )_0$ over the
origin.
By this $T$-action $\rho_k$,
 the natural projection of $E_k$ to $\Bbb A^1$ 
is $T$-equivariant.
Note also that, over $\Bbb A^1$, we have the relative Kodaira embedding 
$$
\mathcal{Z} \;\hookrightarrow\; \Bbb P^*(E_k).
\leqno{(2.2)}
$$
For the structure of $\rho_k$, the following
equivariant trivialization of the vector bundle $E_k$
is known:

\medskip\noindent
{\bf Lemma 2.3} (cf.~\cite{D3}, Lemma 2). 
{\em The holomorphic vector bundle $E_k$ over $\Bbb A^1$ 
can be $T$-equivariantly trivialized by 
$$
E^{}_k \;\cong \; (E^{}_k)^{}_0 \times \Bbb A^1,
$$
where $(E_k)^{}_0$ denotes the fiber of $E_k$ over the origin.
}

\medskip
Let $\lambda_k : T \to \operatorname{GL}_{\Bbb C}((E_k)_0)$ 
denote the algebraic group homomorphism induced 
by $\rho_{k,0}$ on $(E_k)_0$. 
Then the identification in Lemma 2.3 allows us to write
the action $\rho_k$ above in the form
$$
\rho_k(t, (e,s)) \; = \; 
(\lambda_k (t) (e ),\, ts),
\quad (e, s) \in (E_k)^{}_0 \times \Bbb A^1.
\leqno{(2.4)}
$$

\medskip
\section{Proof of Main Theorem}

The isotropy subgroup $\tilde{G}_{\ell}$ of $G_{\ell}$ at 
$[M_\ell ] \in \Bbb P^*(W_{\ell})$ contains
$\hat{G}_{\ell}$ (cf. Section 1) as a subgroup.
Hence
$(M, L^{\ell})$ is Chow-stable if and only if 
$$
|\tilde{G}_{\ell}|<\infty \;\text{ and }\;\text{$G_{\ell}\cdot M_{\ell}$ is closed in $W_{\ell}$,}
$$
because if $\dim \hat{G}_{\ell} < \dim \tilde{G}_{\ell}$, then $\tilde{G}_{\ell}\cdot M_{\ell} =
\Bbb C^* M_{\ell}$, 
and the origin is in the closure of $\tilde{G}_{\ell}\cdot M_{\ell}$
in $W_{\ell}^*$.
Now for all $0< \ell \in \Bbb Z$,  
the identity component $\tilde{G}^0_{\ell}$ of $\tilde{G}_{\ell}$ is
isogenous to an algebraic subgroup  of $\hat{H}$, while
by GIT \cite{Mu}, Proposition 1.5,  
$\tilde{G}^0_m$ is isogenous to $\hat{H}$ for all multiples $m > 0$  
of some fixed integer $\gg 1$.
Hence $(M, L)$ is asymptotically Chow-stable if and only if
$$
\text{$\hat{H}=\{1\}$ and for $\ell \gg 1$, $G_{\ell}\cdot M_{\ell}$ 
is closed in $W^*_{\ell}$.}
\leqno{(3.1)}
$$
 Similarly, if $\ell >0$ is a multiple of some fixed integer $\gg 1$,
we see that the identity component of $\hat{G}_{\ell,k}$ 
with $k \gg 1$ is isogenous to $\hat{H}$.
Hence $(M,L)$ is asymptotically Hilbert-stable if and only if
$$ 
\text{$\hat{H}=\{1\}$ and for all $\ell\gg 1$, 
 $G_{\ell, k}\cdot f_{\ell, k}$ is closed in $F_{\ell, k}$ 
if $k \gg 1$.}
\leqno{(3.2)}
$$
In view of (3.1) and (3.2) above,
(b) of Main Theorem follows immediately from  Fogarty's result
together with (a) of Main Theorem. Hence, we have only to show (a) of Main Theorem.

\medskip
For one-dimensional algebraic torus $T := \Bbb C^*$,
we consider  an algebraic one-parameter 
subgroup 
$$
\lambda :\;  T\; \hookrightarrow \; G_\ell
$$ 
of the reductive algebraic group $G_\ell := \operatorname{SL}(V_\ell)$.
Then to each $\lambda$ as above, we assign a test configuration of exponent 
$\ell$ as follows:

\medskip\noindent
{\em Definition $3.3$}. 
The  {\it DeConcini-Procesi family\/}
(cf. \cite{T1}) associated to $\lambda$ is the 
test configuration of exponent $\ell$ for $(M,L)$
obtained as the $T$-equivariant projective morphism
$$
\pi :\; \mathcal{Z}(\lambda ) \to \Bbb A^1,
$$
where 
$\mathcal{Z}(\lambda)$ is the variety defined as the closure of 
$T\cdot (\iota_\ell (M)\times \{1\})$ in $\Bbb P^*(V_\ell )\times \Bbb A^1$,
and the morphism $\pi$ is induced by the
projection of $\Bbb P^*(V_\ell )\times \Bbb A^1$
to the second factor.
Let $\operatorname{pr}_1: \mathcal{Z}(\lambda ) \to \Bbb P^*(V_\ell )$ 
denote the map induced by the
projection of $\Bbb P^*(V_\ell )\times \Bbb A^1$ to the first factor.
For the open subset $\Bbb C^*\subset\Bbb A^1$,
the holomorphic map $\hbar : \Bbb C^*\to Hilb_{\Bbb P^* (V_{\ell})}$
sending each $t\in \Bbb C^*$ to 
$\hbar (t) :=
\operatorname{pr}_1(\mathcal{Z}(\lambda )_t)\in {Hilb}_{\Bbb P^* (V_{\ell})}$
extends to a holomorphic map 
$$
\tilde{\hbar} : \Bbb A^1 \to Hilb_{\Bbb P^* (V_{\ell})},
$$
where $\mathcal{Z}(\lambda )_s := \pi^{-1}(s)$, $s \in \Bbb A^1$, denotes the scheme-theoretic 
fiber of $\pi$ over $s$.
Now we can regard $\mathcal{Z}(\lambda )$ as the pullback, by $\tilde{\hbar}$, of
the universal family over $Hilb_{\Bbb P^* (V_{\ell})}$.
Note also that $T$ acts on
$\Bbb P^*(V_\ell )\times \Bbb A^1$ by
$$
\;T\times (\Bbb P^*(V_\ell )\times \Bbb A^1) \to (\Bbb P^*(V_\ell )\times 
\Bbb A^1), \quad
(t, (w,s)) \mapsto (\lambda (t) w, t s),
$$ 
where $G_{\ell}$ acts naturally on $\Bbb P^*(V_\ell )$ via the 
contragradient representation.
Then the invertible sheaf 
$$
\mathcal{L} := \operatorname{pr}_1^*\mathcal{O}_{\Bbb P^*(V_\ell )}(1)
$$
over $\mathcal{Z}(\lambda )$ is relatively very ample for the morphism $\pi$, 
and allows us to regard $\pi$ as a projective morphism.
Since the bundle space for $\mathcal{O}_{\Bbb P^*(V_\ell )}(-1)$ 
is identified with the blowing-up of $V_{\ell}^*$ at the origin,
the $G_{\ell}$-action on $V_{\ell}^*$
induces naturally a $T$-action on $\mathcal{L}$ lifting
the $T$-action on $\mathcal{Z} (\lambda )$. 
By restricting $\mathcal{L}$ to 
$\mathcal{Z}(\lambda )_s$, we have isomorphisms
$$
(\mathcal{Z}(\lambda )_s, {\mathcal{L}}_s )
\cong (M, \mathcal{O}_M( L^{\ell})),
\qquad 0 \neq s \in \Bbb A^1,
$$
where $\mathcal{L}_s:= {\mathcal{L}}_{|\mathcal{Z}(\lambda )_s}$ for each $s \in \Bbb A^1$.
Hence $\pi : \mathcal{Z}(\lambda ) \to \Bbb A^1$ is a test configuration of exponent $\ell$ for 
$(M,L)$.

\medskip
For the DeConcini-Procesi family $\pi : \mathcal{Z}(\lambda ) \to \Bbb A^1$ as above,
let $n_k (\lambda ) \in \Bbb Z$ denote the weight of the $T$-action on the 
complex line
$$
\wedge^{m_{k}} (E_k)^{}_0\; 
(\cong \wedge^{m_{k}}H^0(\mathcal{Z}(\lambda )_0, {\mathcal{L}_0^{\,k}})
\;\text{ if $k \gg 1$}),
$$
where $(E_k)^{}_0 := (\pi_*{\mathcal{L}^k})^{}_0$ denotes the fiber, over the origin, 
of the locally free sheaf: 
$\pi_*{\mathcal{L}^k} \to \Bbb A^1$.
If $k \gg 1$, then $\dim H^0(\mathcal{Z}(\lambda )_0, \mathcal{L}_0^{\,k})$ 
is $m_k := \dim H^0(M, \mathcal{O}_M(L^{\ell k}))$,
and we write $m_k$ and $n_k (\lambda )$ as
\begin{align*}
m_k \;\;\;\,\; &= \;\; \Sigma_{i=0}^n\;  \mu_{\ell, i} k^i,
\tag{3.4}\\
n_k(\lambda ) \; &= \;\; \Sigma_{j=0}^{n+1}\; \nu_{\ell, j} (\lambda )k^j,
\tag{3.5}
\end{align*}
where $\mu_{\ell, i}$, $i=0,1,\dots, n$,  and $\nu_{\ell, j} (\lambda )$,
$j=0,1,\dots, n+1$, are rational real numbers
independent of the choice of positive integers $k$. 

\medskip
Let $0 \neq M^{0}_{\ell}\in W_{\ell}^*$ be
such that the associated $[M^{0}_{\ell}]\in \Bbb P^*(W_{\ell})$ 
is the Chow point for the cycle $\mathcal{Z}(\lambda )_0$ 
on $\Bbb P^*(V_{\ell})$ counted with multiplicities.
First, we observe that $\mu_{\ell, n} = \ell^n c_1(L)^n[M]/n! > 0$.
Next, in Section 4, we shall show that
$$
\nu_{\ell, n+1}(\lambda ) \; =\; -\,\frac{a_{\ell}}{(n+1)!}
\leqno{(3.6)}
$$
where $a_{\ell}$ denotes the weight of the $T$-action on $\Bbb C^* M_{\ell}^0$.
We now put
$$
w_{\lambda}(k;\ell) := \; n_k(\lambda ) /(km_k).
$$
{\em Remark}.
Besides the embedding $\mathcal{Z}(\lambda )_0 \hookrightarrow \Bbb P^*(V_{\ell})$,
we also have the embedding
$$
\mathcal{Z}(\lambda )_0 \;\hookrightarrow \; \Bbb P^*((E_1)_0)
$$
for the linear subsystem associated to $(E_1)_0$ in the complete linear system $|\mathcal{L}_0|$ on $\mathcal{Z}(\lambda )_0$.
In the same manner as the weight $a_{\ell}$ above is obtained from the cycle on $\mathcal{Z}(\lambda )_0$ 
on $\Bbb P^*(V_{\ell})$,
we similarly obtain a weight $a'_{\ell}$
from the cycle $\mathcal{Z}(\lambda )_0$ on $\Bbb P^*((E_1)_0)$.
Now by Mumford [6], Proposition 2.11, 
$$
\nu_{\ell, n+1}(\lambda ) = - \frac{a'_{\ell}}{(n+1)!}.
$$
Then (3.6) above claims that $a'_{\ell}$ is replaced by $a_{\ell}$
in this last equality.

\medskip\noindent
{\em Proof of {\rm (a)} of Main Theorem}: 

\medskip
The argument at the beginning of this section shows that the identity 
component $\hat{G}_{\ell_i}^0$ of $\hat{G}_{\ell_i}$ satisfies
$$
\hat{G}_{\ell_i}^0\; \subset \; \tilde{G}^0_{\ell_i},
$$
where $\tilde{G}^0_{\ell_i}$ is isogeneous to an algebraic subgroup of $\hat{H}$. 
Hence the assumption $\hat{H} =\{1\}$ of (a) of Main Theorem implies
$$
|\hat{G}_{\ell_i}| \; <\; \infty\quad
\qquad \text{for all $0 \leq i\in \Bbb Z$}.
$$
Put $K_i:=\Pi_{j=0}^i\, k_j$  
for $0 \leq i \in \Bbb Z$, where we put $K_{-1}:=1$ for simplicity. Moreover, we put 
$\ell_i := \ell K_{i-1}$ for $1 \leq i \in \Bbb Z$. 
Applying Lemma 3.10 below to $(\ell', \ell'', k', k'') = (\ell, \ell_i, K_{i-1}, K_i)$, 
we obtain
$$
w_{\lambda}(K_i; \ell ) > w_{\lambda}(K_{i-1};\ell ),
\qquad i=0,1,2, \dots,
\leqno{(3.7)}
$$ 
for all algebraic one-parameter subgroup 
$\lambda : \Bbb C^* \hookrightarrow G_{\ell}$.
On the other hand, by Appendix, 
we have $n_1 (\lambda ) =0$, i.e., 
$$
w_{\lambda}(K_{-1};\ell )\; =\; w_{\lambda}(1;\ell )\; =\; 0
\leqno{(3.8)}
$$
In view of (3.4) and (3.5), we see that
$$
\lim_{k \to \infty} w_{\lambda}(k; \ell )= \frac{\nu_{\ell, n+1}(\lambda )}{\mu_{\ell, n}}.
\leqno{(3.9)}
$$
By (3.7), (3.8) and (3.9) together with $\mu_{\ell, n} >0$, it follows that
$$
\nu_{\ell, n+1}(\lambda )\; > \; 0
\qquad \text{ for all $\lambda$.}
$$
By (3.6), we conclude that $(M, L^{\ell})$ is Chow-stable,
as required.

\medskip\noindent
{\bf Lemma 3.10}.  {\em Let $n+1 \leq  \hat{k}\in \Bbb Z$,
and let $k'$, $\ell'$ be positive integers with $\ell' \geq n+1$.
 Assume that $G_{\ell''}\cdot f_{\ell'', \hat{k}}$ is closed in $F_{\ell'', \hat{k}}$
for $k'' := \hat{k} k'$ and $\ell'' := k' \ell'$.
If $\hat{H} =\{1\}$, then $w_{\lambda} (k'' ; \ell' )> w_{\lambda}(k'; \ell' ) $ for all 
algebraic one-parameter subgroups $\lambda : \Bbb C^* \hookrightarrow G_{\ell'}$.}

\medskip
\section{Proof of (3.6)}

In this section, we shall prove (3.6) by calculating the term $n_k (\lambda )$ in (3.4) in detail.
Hereafter, by considering
the Decontini-Procesi family
$\mathcal{Z} = \mathcal{Z} (\lambda )$ over $\Bbb A^1$, 
we study the bundles $E_k$, $k =1,2,\dots$ as in Section 2.
A difficulty in calculating $n_k (\lambda )$ comes up when $\mathcal{Z}(\lambda )_0$ 
sits in a hyperplane of $\Bbb P^*(V_{\ell})$. Let 
$N$ be the, possibly trivial,  $T$-invariant maximal linear subspace 
of $V_{\ell}$ vanishing on $\mathcal{Z}(\lambda )_0$, 
where we regard $\mathcal{Z}(\lambda )_0$ as a subscheme
in $\Bbb P^*(V_{\ell}) \; (\cong \{0\}\times \Bbb P^*(V_{\ell}))$. 
Then for some $T$-invariant subspace $Q_1$ of $V_{\ell}$, we write the vector space
$V_{\ell}$ as a direct sum 
$$
V_{\ell} \; =\; Q_1 \oplus N.
$$ 
By $Q_1 = V_{\ell}/N$, we naturally have a $T$-equivariant inclusion
$Q_1 \subset R_1$, where $R_1 := (E_1)_0 =  (\pi_*\mathcal{L})_0\otimes \Bbb C$.
 Then
$$
\mathcal{Z}(\lambda )_0 \subset \Bbb P^*(Q_1) 
\subset \Bbb P^*(V_{\ell}),
$$
i.e., $\mathcal{Z}(\lambda )_0$ sits in the $T$-invariant 
linear subspace $\Bbb P^*(Q_1)$ of
$\Bbb P^*(V_{\ell})$. By taking the direct 
sum of the symmetric 
tensor products for $Q_1$, 
we put $Q := \oplus_{k =0}^{\infty} S^{k}(Q_1)$, 
where 
$S^k(Q_1)$ 
denotes $\Bbb C$ for $k=0$.
Let $J(Q) \subset Q$ denote the $T$-invariant homogeneous ideal of $\mathcal{Z}(\lambda )_0$ 
in $\Bbb P^*(Q_1)$.
Then by setting $J(Q)_k := J(Q) \cap S^k(Q_1)$, 
we define 
$$
Q_k := S^k(Q_1)/J(Q)_k.
$$
By Theorem 3 in \cite{Mum}, the natural homomorphism: 
$S^k(E_1) \to E_k$ is surjective over $\Bbb A^1 \setminus \{0\}$
for all positive integers $k$. We also have the $T$-equivariant inclusion
$$
Q_k \; \subset \; R_k,
\qquad 0 < k\in \Bbb Z,
\leqno{(4.1)}
$$
where $R_k := (E_k)_0\; =\; (\pi_*\mathcal{L}^k)_0\otimes\Bbb C$.
By choosing general elements $\sigma_i$, $i =0,1,2,\dots n$, in $Q_1$,
we have a surjective holomorphic map
$$
\operatorname{pr}^{}
_{\Bbb P^n}:\,\mathcal{Z}(\lambda )_0 \to \Bbb P^n(\Bbb C),
\qquad  z \mapsto 
(\sigma_0(z):\sigma_1(z): \dots :\sigma_n (z)),
$$
so that the fiber $\operatorname{pr}_{\Bbb P^n}^{\,-1}(q)$ over $q:= (1:0:0 \dots :0)$
consists of
$r$ points counted with multiplicities, where $r:= \ell^n c_1(L)^n [M]$. 
For each $k \gg 1$, we
consider the subspace $F_k:= \operatorname{pr}^*H^0(\Bbb P^n (\Bbb C ), \mathcal{O}^{}_{\Bbb P^n}(k))$
of $Q_k$. Then
$$
\dim \; F_k \; =\; \frac{(n+k)!}{n!\,k!}
$$
is a polynomial in $k$ of degree $n$ with leading coefficient $1/n!$.
For some positive integer $k_0$,  there exist elements
$\tau_1$, $\tau_2$, \dots , $\tau_r$ in $Q_{k_0}\setminus F_{k_0}$  
which separate the points in $\operatorname{pr}_{\Bbb P^n}^{\,-1}(q)$ 
including infinitely near points. Then for $k \gg 1$, the linear subspaces
$$
\tau_1 F_{k-k_0}, \; \tau_2 F_{k-k_0}, \dots , \tau_r F_{k-k_0}
$$
of $Q_k$ altogether span a linear subspace of dimension
$$
r \dim\; F_{k-k_0} \; =\; r\, \frac{(n+k - k_0)!}{n!\, (k-k_0)!} \; 
=\; \frac{r}{n!} \,k^n + \text{ lower order term in $k$}.
$$     
In view of (4.1),  
$r\dim\, F_{k-k_0} \leq\dim\, Q_k \leq \dim\, R_k = m_k$. Hence
$$
\dim\; R_k/Q_k \; \leq \; C_1 \,k^{n-1}
\leqno{(4.2)}
$$
for some positive constant $C_1$ independent of $k$. 
Put $\delta_k := \dim Q_k$, and 
let $q_k (\lambda )$ denote the weight of the $T$-action on $\wedge^{\delta_k} Q_k$,
where the weight of the $T$-action on $\wedge^{m_k} R_k$ is $n_k(\lambda )$.
Then the weight of the $T$-action on $\wedge^{m_k -\delta_k} (R_k/Q_k)$ 
is $n_k (\lambda ) - q_k (\lambda )$.
On the other hand, in view of Remark 4.6 below, 
the weight $\alpha$ for the $T$-action on 
every $1$-dimensional $T$-invariant subspace $A$
of $R_k/Q_k$ satisfies
$$
|\alpha | \; \leq \; C_2 \, k
\leqno{(4.3)}
$$
for some positive constant $C_2$ independent of the choice of $k$.
Then we see from (4.2) and (4.3) that
$$
|n_k (\lambda ) - q_k (\lambda ) |\; \leq \; C_1 C_2 \, k^n.
\leqno{(4.4)}
$$
Now a classical result of Mumford \cite{Mu0}, Proposition 2.11, asserts that 
$$
q_k (\lambda )\;  =\; -\frac{a_{\ell}\, k^{n+1}}{(n+1)!} + \text{ lower order term 
in k},\leqno{(4.5)}
$$
where the weight in \cite{Mu0} and ours have opposite sign.
From (3.5), (4.4) and (4.5), we obtain (3.6) as required.

\medskip\noindent
{\em Remark $4.6$}.
Put $X_0:= \mathcal{Z}(\lambda )_0$. For $X_0$ sitting in 
$\Bbb P^* (V_{\ell})\times \{0\}$, we choose a sequence of scheme-theoretic 
intersections
$$
X_j := X_0 \cap \Sigma_1 \cap \Sigma_2 \cap \dots \cap \Sigma_j,
\qquad j= 1,2,\dots,n,
$$ 
where $\Sigma_1$, $\Sigma_2$, \dots, $\Sigma_n$ are $n$ distinct general 
hyperplanes in $\Bbb P^* (V_{\ell})\times \{0\}$.  
Then there exists an integer $i_0$
satisfying $i_0 \gg n$ such that
$$
H^p(X_j, \mathcal{O}_{X_j}(\mathcal{L}_0^{\,i})) = 
\{0\},
\qquad  i \geq i_0 -n, 
$$
for all $p >0 $ and $j = 0,1,\dots,n$. Then by the arguments in the 
proof of Theorem 2 in \cite{Mum}, the natural homomorphisms
$$
H^0(\mathcal{Z}(\lambda )_0, 
\mathcal{L}^{\,i}_0)\otimes H^0(\mathcal{Z}(\lambda )_0, \mathcal{L}_0)
\; \to\; H^0(\mathcal{Z}(\lambda )_0, 
\mathcal{L}^{\,i+1}_0),
\quad i\geq i_0,
$$
are surjective. In particular, for all positive integers $k$, 
the natural homomorphisms: $(R_{i})^{\otimes k} 
\to R_{ik}$ are surjective for all integers $i \geq  i_0$.

\medskip
\section{Proof of Lemma 3.10}

In this section,  we
apply Section 2 to $\mathcal{Z} = \mathcal{Z}(\lambda )$ 
and $\ell = \ell'$, where the actions of $\,T := \Bbb C^*$ on $\mathcal{L}$
and $\mathcal{Z}(\lambda )_0$
are induced by the one-parameter group $\lambda : \Bbb C^* \hookrightarrow G_{\ell'}$ 
in Lemma 3.10, where for each positive integer $k$, the corresponding $T$-action $\rho_k$  
on $E_{k}$ induces
$$
\lambda_k\,:\, T \to \operatorname{GL}_{\Bbb C}((E_k)_0)
$$
as in (2.4).
Recall that $k'' = \hat{k} k'$ and $\hat{k} \geq n+1$.  
For each $s \in \Bbb A^1\setminus\{0\}$, 
let $\,\mathcal{I}(\mathcal{Z}(\lambda )_s )\,$ denote
the kernel  
of the naural $T$-equivariant surjective homomorphism
$$
S^{\hat{k}}(E_{k'})_s \to (E_{k''})_s,
$$ 
between fibers over $s$ for bundles $S^{\hat{k}}(E_{k'})$ and $E_{k''}$,
where the $T$-actions on $S^{\hat{k}}(E_{k'})$ and $E_{k''}$ 
are by $\rho_{k'}$ and $\rho_{k''}$, respectively.
By the trivialization in Lemma 2.3 applied to $k = k'$, we can identify each $\Bbb G\operatorname{r}_s$, $s \in \Bbb A^1$, 
with $\Bbb G\operatorname{r}_0$. Here $\Bbb G\operatorname{r}_s$ denotes the 
complex Grassmannian  
of all complex $\gamma_{\ell'',\hat{k}}$-planes through the origin in $S^{\hat{k}}(E_{k'})_s$.
The holomorphic map sending each $s \in \Bbb A^1 \setminus \{0\}$ to 
$\,\mathcal{I}(\mathcal{Z}(\lambda )_s )\,$ regarded as an element in 
($\Bbb G\operatorname{r}_s \cong ) \, \Bbb G\operatorname{r}_0$ extends naturally to 
a holomorphic map: $\Bbb A^1 \to \Bbb G\operatorname{r}_0$, where the image 
of the origin under this holomorphic map will be denoted
by $\,\mathcal{I}(\mathcal{Z}(\lambda )_0 )\,$ by abuse of terminology.
For the inclusion
$\mathcal{Z}(\lambda ) \hookrightarrow\Bbb P^* (E_{k'})$ in (2.2),
the action of each $t\in T$ maps
$\mathcal{Z}(\lambda )_{s}$ onto $\mathcal{Z}(\lambda )_{ts}$, 
and we have 
$$
\mathcal{I}(\mathcal{Z}(\lambda )_{ts}) \; =\; 
\rho_{k'} (t) ( \mathcal{I}(\mathcal{Z}(\lambda )_{s})).
\leqno{(5.1)}
$$
Here, via the $T$-action on $\mathcal{Z}(\lambda )_0$, 
$T$ acts on $S^{\hat{k}}(E_{k'})_0$
preserving $\mathcal{I}(\mathcal{Z}(\lambda )_0)$.
At $s =1$, the fiber $\mathcal{Z}(\lambda )_s := \pi^{-1}(s)$ over $s$
is thought of as
$\iota_{\ell'}(M)$ sitting in $\Bbb P^*(V_{\ell'})$. 
Hence by the notation in Section 1,  we have 
$$
\mathcal{I}(\mathcal{Z}(\lambda )_s )_{|s=1}\; =\; I_{\ell'' , \hat{k}}
$$ 
by identifying
${E_{k'}}_{|s=1}$, ${E_{k''}}_{|s=1}$ with $V_{\ell''}$, $V_{\ell'' \hat{k}}$, 
respectively. 
Consider the closed disc $\Delta := \{ s\in \Bbb A^1\,;\, |s|\leq 1 \,\}$
of $\Bbb A^1$.
Since
$\wedge^{\gamma_{\ell'' , \hat{k}}} \mathcal{I}(\mathcal{Z}(\lambda )_s )$ in 
$\wedge^{\gamma_{\ell'' , \hat{k}}} S^{\hat{k}}(E_{k'})_s$ is a complex line, 
for each $s \in \Delta$, 
we can choose 
an element $\psi_{\ell'',\hat{k}}(s) \neq 0$ in the line in such a way that
$\psi_{\ell'',\hat{k}}(s) $ depends on
$s$ holomorphically.
Then $\psi_{\ell'',\hat{k}}(1)$ is regarded as a nonzero element in the line 
$\wedge^{\gamma_{\ell'' , \hat{k}}} I_{\ell'' , \hat{k}}$ in 
$\wedge^{\gamma_{\ell'' , \hat{k}}} S^{\hat{k}} (V_{\ell''})$.
By the trivialization in Lemma 2.3 applied to $k = k'$, 
we hereafter identify $(E_{k'})_s$,
$s \in \Delta$,  
with $(E_{k'})_0$. Consequently, this identification allows us to regard
$\psi_{\ell'',\hat{k}}(s)$ as an element 
in $\Psi := \wedge^{\gamma_{\ell'' , \hat{k}}} S^{\hat{k}}(E_{k'})_0$ for each 
$s \in \Delta$, and
$G_{\ell''}$ is viewed as $\operatorname{SL}_{\Bbb C}((E_{k'})_0)$.
For each $t$, $t'\in \Bbb C^*$, by taking an unramified
cover of $\Bbb C^*$ of degree $m_{k'}$, we can write
$$
t = \tilde{t}_{}^{m_{k'}}\quad
\text{ and }\quad
t' = (\tilde{t}')_{}^{m_{k'}},
$$
for $\tilde{t}$, $\tilde{t}' \in \Bbb C^*$, where $m_{k'}$ is the rank of the
vector bundle $E_{k'}$.
The closedness of $G_{\ell''}\cdot f_{\ell'',\hat{k}}$ 
in $F_{\ell'',\hat{k}}$ in the assumption of Lemma 3.10 means that 
the orbit $\operatorname{SL}_{\Bbb C}((E_{k'})_0)\cdot 
\psi_{\ell'',\hat{k}}(1)$ is closed in $\Psi$.
Now by the Hilbert-Mumford stability criterion, 
$$
\text{ $\hat{\lambda} (\Bbb C^*)\, \cdot  \psi_{\ell'',\hat{k}}(1)$ is closed in $\Psi$,}
\leqno{(5.2)}
$$
where $\hat{\lambda} : \Bbb C^*\to \operatorname{SL}_{\Bbb C}((E_{k'})_0)$ is an algebraic group 
homomorphism defined by
$$
\hat{\lambda} ( \tilde{t} ) :=\; \frac{\lambda_{k'} (t)\;}{\det \lambda_{k'} (\tilde{t})},
\qquad \tilde{t}\in \Bbb C^*, 
$$
for $\lambda_{k'}$ as in Section 2.
To each $\psi_{\ell'',\hat{k}}(s)$,  $s \in \Delta$, we can naturally assign an 
element $[\psi_{\ell'',\hat{k}}(s)]$ in the complex Grassmannian 
$\Bbb G\operatorname{r}_0$.
Here $[\psi_{\ell'',\hat{k}}(s)]$ corresponds to the subspace
$\mathcal{I}(\mathcal{Z}(\lambda )_s)$ in $S^{\hat{k}}(E_{k'})_0$
via the identification $S^{\hat{k}}(E_{k'})_s\cong S^{\hat{k}}(E_{k'})_0$
in terms of the trivialization in Lemma 2.3 applied to $k = k'$.
Obviously,
$$
\text{$[\psi_{\ell'',\hat{k}}(s)] \,\to\, [\psi_{\ell'',\hat{k}}(0)]$\; as $s \to 0$.}
$$
Moreover, in view of (5.1), we obtain
$$
\lambda_{k'}(t)\cdot \psi_{\ell'',\hat{k}}(s)\; \in \;\Bbb C^* 
\cdot \psi_{\ell'',\hat{k}}(ts), 
\qquad s \in \Delta,
$$
for all $t\in \Bbb C^*$ satisfying $|t|\leq 1$.   For some $\varepsilon \in \Bbb R$
with $0 <\varepsilon \ll 1$, we put 
$D_{\varepsilon} := 
\{\,t\in \Bbb C^*\,;\, |t| <\varepsilon \,\}$. Then  
$$
\hat{\lambda} (\tilde{t})\cdot \psi_{\ell'',\hat{k}}(1)\; =\;
\frac{\lambda_{k'}(t)\cdot \psi_{\ell'',\hat{k}}(1)}{\det \lambda_{k'}(\tilde{t})}\;
=\; \tilde{t}_{}^{\beta}\psi ( t ),
\qquad  t \in 
D_{\varepsilon},
\leqno{(5.3)}
$$
for some nonvanishing holomorphic map $\Delta_{\varepsilon} \owns s \mapsto \psi (s) \in \Psi$, where by
$\Delta_{\varepsilon}$, we mean the subset $\{\, s \in \Bbb C\,;
\, |s| \leq \varepsilon\,\} = D_{\varepsilon} \cup \{0\}$ of $\Delta$. 
Now by (5.2) and (5.3), we obtain
$$
\beta \; < \; 0.
\leqno{(5.4)}
$$  
On the other hand, since the map $\psi$ is continuous, (5.3) implies
$$
\lim_{t\to 0}\;\tilde{t}^{-\beta}\,
\hat{\lambda} (\tilde{t})\cdot \psi_{\ell'',\hat{k}}(1)\; =\; \psi (0).
\leqno{(5.5)}
$$
If $t$, $t' \in D_{\varepsilon}$, then from (5.3),
it follows that
$$
\hat{\lambda} (\tilde{t})\hat{\lambda} (\tilde{t}') \cdot \psi_{\ell'',\hat{k}}(1) \; =\; 
\hat{\lambda} (\tilde{t}\tilde{t}')\cdot \psi_{\ell'',\hat{k}}(1)\; 
=\;  (\tilde{t}\tilde{t}')^{\beta} \psi (tt').
$$
Hence $\hat{\lambda} (\tilde{t})\{\tilde{t'}_{}^{-\beta}
\hat{\lambda} (\tilde{t}') \cdot \psi_{\ell'',\hat{k}}(1)\}
= \tilde{t}^{\beta} \psi (tt')$. Let $ t' \to 0$. Then this together with (5.5)  implies
$$
\hat{\lambda} (\tilde{t}) \cdot \psi (0)
\; =\; \tilde{t}^{\beta}\psi (0).
\leqno{(5.6)}
$$
In view of (5.4) and (5.6),  the argument 
as in \cite{D1}, 2.3, applied to $S^{\hat{k}}(E_{k'})_0 \to (E_{k''})_0$
allows us to obtain
$$
0\; >\; \beta\;
=\; \frac{\hat{k} m_{k''} n_{k'}(\hat{\lambda} )}{m_{k'}} - n_{k''}(\hat{\lambda} )\; 
=\; k'' m_{k''}\left\{\, \frac{n_{k'}(\hat{\lambda})}{k'm_{k'}} - 
\frac{n_{k''}(\hat{\lambda})}{k''m_{k''}}\,\right\},
$$
where $n_{k'}(\hat{\lambda})$ and $n_{k''}(\hat{\lambda})$ are 
the weights of the 
$\Bbb C^*$-actions on $\wedge^{m_{k'}}(E_{k'})_0$ and 
$\wedge^{m_{k''}}(E_{k''})_0$, respectively, induced by $\hat{\lambda}$.
Since $\lambda_{k'}$ is induced by $\lambda$, the definition of $\tilde{t}$ and 
$\hat{\lambda}$ 
shows that
$$
\frac{n_{k'}(\hat{\lambda})}{k'm_{k'}} - 
\frac{n_{k''}(\hat{\lambda})}{k''m_{k''}}
\; =\; m_{k'}\{\, w_{\lambda}(k'; \ell' ) - w_{\lambda}(k'' ; \ell' )\,\}
$$
and hence 
$w_{\lambda}(k'; \ell' ) < w_{\lambda}(k'' ; \ell' )$,
as required. This completes the proof of Lemma 3.10.

\medskip
\section{Appendix}

The purpose of this section is to study the $T$-action $\rho_{k,0}$ on $(E_k)_0$ 
with $k=1$ for 
the DeConcini-Procesi family 
$\mathcal{Z} = \mathcal{Z}(\lambda )$ over $\Bbb A^1$.
Let
$$ 
\text{$\tilde{\operatorname{pr}}_1: \Bbb P^* (V_{\ell} ) \times \Bbb A^1 
\to \Bbb P^* (V_{\ell} )$,
\quad $\tilde{\pi} : \Bbb P^* (V_{\ell} ) \times \Bbb A^1 \to \Bbb A^1$}
$$ 
be the projections to respective factors.
Put $\tilde{\mathcal{L}}:= 
\tilde{\operatorname{pr}}_1^*\mathcal{O}_{\Bbb P^*(V_{\ell})}(1)$.
Then for every $e\in V_{\ell}$,
the map sending each $s \in \Bbb A^1$ to 
$(e, s) \in V_{\ell} \times \Bbb A^1$ defines a holomorphic section, denoted by $\tau (e)$, 
 in
$H^0(\Bbb A^1, \tilde{\pi}_*\tilde{\mathcal{L}})$.
The pullback $\iota^*\tau (e)$ by the inclusion map
$$
\iota : \; \mathcal{Z}(\lambda ) \; \hookrightarrow \Bbb P^* (V_{\ell} ) \times \Bbb A^1
$$
is naturally regarded as a holomorphic section of $E_1$ over $\Bbb A^1 = \{s \in \Bbb C\}$,
where $s$ is the affine coordinate for $\Bbb A^1$.
Note that, for $s \neq 0$, the map
$$
V_{\ell} \owns e \; \mapsto \; \{\iota^*\tau (e)\}(s) \in (E_1)_s
$$
is a linear isomorphism. Here $E_k$ with $k =1$ is written as $E_1$,
and $(E_1)_s$ denotes the fiber of the vector bundle $E_1$ over $s$.
In terms of the $T$-action on $V_{\ell}$ via
the one-parameter group $\lambda : T \to \operatorname{SL}(V_{\ell})$, 
write the vector space $V_{\ell}$ as a direct sum
$$ 
N \; =\; \bigoplus_{i=1}^p \; N_i,
\leqno{(6.1)}
$$
where $N_i = \{\, e\in V_{\ell}; \lambda (t) e = t^{\alpha_i}e\,\text{ for all $t\in T$}\}$ 
 for some mutually distinct integers $\alpha_i$,  $i =1,2,\dots, p$.
For each $i$, consider the $\Bbb C [s]$-module $N_i [s] := N_i \otimes_{\Bbb C} \Bbb C [s]$,
where by $\Bbb C [s]$, we mean the ring of polynomials in $s$ with coefficients in $\Bbb C$.
Let $\{ e_1, e_2, \dots, e_{n_i}\}$ be a basis for the vector space $N_i$,
where $n_i := \dim N_i$.
For every $e\in N_i[s]$, by writing $e$ as a sum
$\Sigma^{n_i}_{j=1}f_j(s) e_j\in N_i [s]$ for some polynomials $f_j(s) \in \Bbb C [s]$ 
in $s$, we put
$$
\tau (e) := \sum_{j=1}^{n_i} f_j (s) \tau (e_j) \in 
H^0(\Bbb A^1, \tilde{\pi}_*\tilde{\mathcal{L}}).
$$
From the $T$-action on $V_{\ell}$ via $\lambda$,
we have a natural fiberwise 
$T$-action 
on the trivial bundle $V_{\ell} \times \Bbb A^1$ over $\Bbb A^1$.
This then induces a fiberwise $T$-action on the vector bundle $E_1$ over $\Bbb A^1$,
while the restriction of this induced $T$-action to the fiber $(E_1)_0$ is 
exactly $\rho_{1,0}$.
 
\medskip\noindent
{\bf Proposition 6.2.} 
{\em There exists a non-decreasing sequence of nonnegative integers
$\beta_{i1} \leq \beta_{i2} \leq \dots \leq \beta_{in_i}$ 
together with
$\Bbb C [s]$-generators $\{\, e_{ij}\,;\, j=1,2, \dots, n_i\}$ 
for the $\Bbb C [s]$-module $N_i [s]$ such that
$$
\iota^* \tau (e_{ij})\;  =\; s^{\beta_{ij}} \sigma_{ij}, 
\qquad i=1,2,\dots,p,\; j=1,2,\dots, n_i,
\leqno{(6.3)}
$$
for some holomorphic sections 
$\sigma_{ij}$ to $E_1$ over $\Bbb A^1$, where 
$$
\{\,\sigma_{ij}(0)\,;\, i = 1,2,\dots, p,\,j=1,2,\dots, n_i\,\}
\leqno{(6.4)}
$$
forms a basis for the vector space $(E_1)_0$.
}

\medskip\noindent
{\em Proof}: 
By induction on $j = 1,2, \dots, n_i$, we 
define $e_{ij}$ and $\sigma_{ij}$ from
$\{\, e_{i1}, e_{i2}, \dots, e_{i\,j-1}\,\}$ and 
$\{\, \sigma_{i1}, \sigma_{i2}, \dots, \sigma_{i\,j-1}\,\}$
as follows. Let $B_{j-1}$ be the $\Bbb C [s]$-submodule of $N_i[s]$ 
generated by $\{\, e_{i1}, e_{i2}, \dots, e_{i\,j-1}\,\}$, where 
we put $B_{j-1}=\{0\}$ for $j =1$.
Let $\mathcal{Y}_{ij}$ denote the set of all $\Bbb C [s]$-submodules $Y\subset N_i [s]$ 
generated by $n_i - j +1$ elements such that 
$$
Y + B_{j-1} = N_i [s],
\leqno{(6.5)}
$$
where $Y + B_{j-1}$ is the $\Bbb C [s]$-submodule of $N_i [s]$ generated by 
$Y$ and $B_{j-1}$.
For each $Y \in \mathcal{Y}_{ij}$, let $\beta (Y)$ denote the maximal nonnegative integer 
$\beta$ such that all $\iota^*\tau (e)$, $e\in Y$, are divisible by $s^{\beta}$ in 
$H^0(\Bbb A^1, \mathcal{O}_{\Bbb A^1}(E_1))$. 
In view of the inequality $j \leq n_i$, the maximum
$$
\beta_{ij} := \max_{Y\in 
\mathcal{Y}_{ij}} \beta (Y)
$$ 
exists because, otherwise, (6.5) would imply that
$\iota^*\tau (N_i)  \subset \iota^*\tau (B_{j-1})$ modulo $s^{\beta}$ 
for all positive integers $\beta$, in contradiction to
$n_i > j-1$. By the definition of $\beta_{ij}$, it now easily follows that 
$\beta_{i1} \leq \beta_{i2} \leq \dots \leq \beta_{in_i}$.
Take an element $Y_{ij}$ of $\mathcal{Y}_{ij}$ 
such that $\beta (Y_{ij}) =\beta_{ij}$.
Then the maximality of $\beta_{ij}$ allows us to obtain $e_{ij}\in Y_{ij}$ 
and $\sigma_{ij} \in H^0(\Bbb A^1, \mathcal{O}_{\Bbb A^1}(E_1))$ satisfying
$\iota^*\tau (e_{ij}) = s^{\beta_{ij}} \sigma_{ij}$ such that $\sigma_{ij}(0)$ is
$\Bbb C$-linearly independent from 
$\sigma_{i1}(0)$, $\sigma_{i2}(0)$, \dots, $\sigma_{i\,j-1}(0)$ in $(E_1)_0$.
Since this induction procedure stops at $j = n_i$, we obtain
both (6.3) and the required condition for (6.4).

\medskip
Now the vector bundle $E_1$ is generated by 
the global sections $\{\,\sigma_{ij}\,;\, i = 1,2,\dots, p,\,j=1,2,\dots, n_i \,\}$
 over $\Bbb A^1$. 
Then by (6.1) and (6.3), 
$$
\rho_{1,0} (t, \,\sigma_{ij}(0) ) \; =\;
t^{\alpha_i} \sigma_{ij}(0)
$$
for all $i$ and $j$.
In particular, $n_1(\lambda ) = \Sigma_{i=1}^p n_i\alpha_i$.
Since $\lambda$ is an algebraic one-parameter subgroup in $G_{\ell} = 
\operatorname{SL}(V_{\ell})$, 
by the definition of $N_i$, it follows from (6.1) that 
$$
1 = \det (\lambda (t)) = t^{\Sigma_{i=1}^p n_i\alpha_i}
$$ 
for all $t \in T$, i.e., $n_1 (\lambda ) = 0$.
Note that this equality follows also from Lemma 2.3 by the equivariant 
isomorphism $(E_1)_1 \cong (E_1)_0$ (see also \cite{O}).


\bigskip\noindent
{\footnotesize
{\sc Department of Mathematics}\newline
{\sc Osaka University} \newline
{\sc Toyonaka, Osaka, 560-0043}\newline
{\sc Japan}}
\end{document}